\documentclass[12pt]{article}
\usepackage{amsmath,amssymb,amsthm}

\numberwithin{equation}{section}

\newtheorem{thm}{Theorem}[section]
\newtheorem{lem}{Lemma}[section]

\newcommand{\n}{\nonumber}

\newcommand{\si}{\sigma_R }

\renewcommand{\a}{\alpha}

\newcommand{\vare}{\varepsilon}
\newcommand{\s}{\sigma}

\renewcommand{\O}{\Omega}
\newcommand{\bb}{\begin{equation}}
\newcommand{\ee}{\end{equation}}
\newcommand{\bq}{\begin{eqnarray}}
\newcommand{\eq}{\end{eqnarray}}
\newcommand{\bqn}{\begin{eqnarray*}}
\newcommand{\eqn}{\end{eqnarray*}}
\newcommand{\R}{\mathbb{R}}
\newcommand{\I}{\infty}
\newcommand{\pd}{\partial}

\newcommand{\sz}{S_0} 

\begin{document}
\title{Remarks on the asymptotically discretely self-similar solutions of the Navier-Stokes\\ and the Euler equations}
\author{Dongho Chae   \\
Chung-Ang University\\ Department of Mathematics\\
 Seoul 156-756, Republic of Korea\\
e-mail: dchae@cau.ac.kr}
\date{}
\maketitle

\begin{abstract}
We study scenarios of self-similar type blow-up for the incompressible Navier-Stokes and the Euler equations.
The previous notions of the discretely (backward) self-similar solution and the asymptotically self-similar solution are generalized to the locally asymptotically discretely self-similar solution.
We prove that there exists no such locally asymptotically discretely self-similar blow-up for the 3D Navier-Stokes equations if the blow-up profile is a time periodic function belonging to $C^1(\Bbb R ; L^3(\Bbb R^3)\cap C^2 (\Bbb R^3))$.
For the 3D Euler equations we show that the scenario of asymptotically discretely self-similar blow-up  is excluded  if the blow-up profile satisfies suitable integrability conditions.
\end{abstract}
\section{The main theorems}
\setcounter{equation}{0}
\subsection{Navier-Stokes equations}
We consider  the Cauchy problem of  the 3D Navier-stokes equations.
$$
(NS) \left\{\aligned  &\partial_t v+v\cdot \nabla v =-\nabla p+\Delta v,\label{nv1}\\
&\mathrm{div }\, v=0,\\
& v(x,0)=v_0(x),
\endaligned \right.
$$
 where $v(x,t)=(v_1 (x,t), v_2 (x,t), v_3 (x,t))$ is the velocity, $p=p(x,t)$ is the
  pressure, and $v_0(x)$ is the initial data satisfying div$v_0 =0$.
  We study the possibility of finite time blow-up of smooth solution of (NS). By translation of time we may assume that
  the solution is smooth for $t<0$, and the blow-up happens at $t=0$.
 We say a solution  $v(x,t)$ to (NS) is a (backward) self-similar blowing up solution at $t=0$ if there exists $(V,P)$ such that
 \bb\label{ss}
v(x,t) =\frac{1}{\sqrt{-t}} V\left(\frac{x}{\sqrt{-t}}\right),\quad p(x,t)= -\frac{1}{t} P\left(\frac{x}{\sqrt{-t}}\right),
 \ee
  For such solution we have  $\lambda v(\lambda x,\lambda^2 t)=v(x,t)$ for all $\lambda \in \Bbb R$ and for all $(x,t)\in \Bbb R^3\times (-\infty, 0)$,
 The nonexistence of nontrivial solution  given by (\ref{ss}) was shown by Ne\v{c}as-R{\aa}\v{z}i\v{c}ka-\v{S}ver\'{a}k(\cite{nrs}), and Tsai(\cite{tsai})(see also \cite{mil}). It can also be deduced from the result of \cite{esc}.
 These results are generalized by introduction of the more general notion of the asymptotically self-similar blow-up solutions, and their exclusion in \cite{cha1}, which was motivated by earlier study of asymptotically self-similar solutions for the nonlinear heat equation by Giga and Kohn\cite{gig1, gig2}(see also \cite{hou}).
  A different direction of generalization of the notion can be done  by considering  the discretely (backward) self-similar blow-up at $(x,t)=(0,0)$ for a solution of (NS) is a solution $v(x,t)$, which means that  there exists  $\lambda \neq 1$ such that  $\lambda v(\lambda x, \lambda ^2 t)=v(x,t)$ for all  $(x,t)\in \Bbb R^3 \times (-\infty, 0)$. We note that the notion of forward discretely self-similar solutions to (NS) was studied earlier by Tsai in \cite{tsai0}(see also \cite{egg}).
If we make a self-similar transform, $(v,p)\to  (V,P)\in [C^\infty (\Bbb R^{3+1})]^2$ by the formula,
\bb\label{nsds}
v(x,t) =\frac{1}{\sqrt{-t}} V(y,s),\quad p(x,t)= -\frac{1}{t} P(y,s),
 \ee
 where
$ y=x/ \sqrt{-t}$ and $ s=-\log (-t)$,
then, $(V,P)$ satisfies
\begin{equation}\label{leray} \left\{
\aligned
&V_s+ \frac {1}{2} V +\frac {1}{2}(y \cdot \nabla)V + (V\cdot \nabla )V =-\nabla
P +\Delta V,\\
& \mathrm{div}\,V=0,
\endaligned \right.
\end{equation}
which has the periodicity in time,
$V(y,s)= V(y, s+S_0)$ for all $(y,s)\in \Bbb R^{3+1}$ with $S_0= -2\log \lambda.$
We first prove the nonexistence of discretely  self-similar blow-up for the $L^3(\Bbb R^3)$ solution to (NS), which is an easy consequence of \cite{ser}, as is shown in the next section.
\begin{thm}
 If $v\in C(-\infty, 0; L^3 (\Bbb R^3))$ a solution to (NS), which  blows up at $t=0$, then $t=0$ is not a time for discretely self-similar blow up.
\end{thm}

Next we consider more general possibility of locally asymptotically discretely self-similar blow-up.
 Let $q\in [1, \infty)$, and $v\in C (-\infty, 0; L^q_{loc}(\Bbb R^3))$ be  a solution of (NS), which blows up at $t=0$.  We say that $(x,t)=(0,0)$ is  a space-time point of {\em locally asymptotically
 discretely self-similar blow-up in the sense of $L^q$} if there exist $R>0$ and  a solenoidal vector field $\overline{V}(y,s)\in C(\Bbb R; L^q (\Bbb R^3))$ with $\overline{V}(y,s)=\overline{V}(y,s+S_0)$ for some $S_0\neq 0$ and for all $(y,s)\in \Bbb R^3\times \Bbb R$  such that
\bb\label{conv1}
\lim_{t\uparrow 0} (-t)^{\frac{q-3}{2q}}  \sup_{t<\tau<0} \left\| v(\cdot ,\tau) -\frac{1}{\sqrt{-\tau}} \overline{V}\left(\frac{\cdot }{\sqrt{-\tau}}, -\log (-\tau)\right)\right\|_{L^q (B(0, R \sqrt{-t}))}=0.
\ee
Note that this is much more general notion than the discrete self-similarity, which corresponds to the equality inside of the norm in (\ref{conv1}) for all $R>0$.
The following is a generalization of Theorem 1.2 of \cite{cha1} and Theorem 1.1 above.
\begin{thm}
 Let $v\in C (-\infty, 0; L^3_{loc}(\Bbb R^3))$  be a solution to (NS), which blows up at $t=0$, then it is not a time for locally asymptotically discretely self-similar blow-up in the sense of $L^q(\Bbb R^3)$ if  $q
 \in [2, \infty]$, and the blow-up profile satisfies $ \overline{V}\in C (\Bbb R; L^3 (\Bbb R^3)\cap C^2 (\Bbb R^3))$.
\end{thm}
\noindent{\em Remark 1.1 } Although we considered possible blow-up at $(x,t)=(0,0)$ in the above, the result holds also for any  $(x,t)$  by translation.\\

\subsection{The  Euler equations}
The question of the finite time blow-up for the 3D Euler equations is also one of the most outstanding
problem in the mathematical fluid mechanics(see \cite{maj} for the physical and mathematical backgrounds on the problem and \cite{bea, con} for well-known partial results  in this direction).
Here we study the possibility of existence of self-similar type blow-up for a solution to the 3D Euler equations.
$$
(E) \left\{\aligned  &\partial_t v+v\cdot \nabla v =-\nabla p,\label{e1}\\
&\mathrm{div }\, v=0,\\
& v(x,0)=v_0(x)
\endaligned \right.
$$
We say a solution $v$ is a self-similar blowing up solution at $t=0$ if there exists $(V, P)$ such that
\bb\label{ss} v(x,t)= \frac{1}{(-t)^{\frac{\alpha}{\alpha+1}}} V\left(\frac{x}{(-t)^{\frac{1}{\alpha+1}}}\right),\quad
p(x,t)= \frac{1}{(-t)^{\frac{2\alpha}{\alpha+1}}} P\left(\frac{x}{(-t)^{\frac{1}{\alpha+1}}}\right).
\ee
For physical motivation on this type scenarios we refer \cite{pel}.
For such $(v,p)$ we have the following scaling invariance $v(x,t)= \lambda ^\alpha v(\lambda x, \lambda ^{\alpha+1} t)$ and
$p(x,t)= \lambda ^{2\alpha} v(\lambda x, \lambda ^{\alpha+1} t)$ for all  $\lambda >0$, $\alpha \in \Bbb R$, and for all $(x,t)\in \Bbb R^3\times (-\infty, 0)$. The nonexistence of nontrivial self-similar blowing up solution under suitable assumption on the blow-up profile $V$ was obtained in \cite{cha3}.
A discretely self-similar solution $v$ is a solenoidal vector field, for which there exist $\lambda >0, \alpha >-1$ such that
           $\label{cond1} v(x,t)= \lambda ^\alpha v(\lambda x, \lambda ^{\alpha+1} t)\quad \forall (x,t)\in \Bbb R^3 \times (-\infty, 0)$.
If we represent $(v,p)$ by
\bb\label{ss} v(x,t)= \frac{1}{(-t)^{\frac{\alpha}{\alpha+1}}} V(y,s),\quad
p(x,t)= \frac{1}{(-t)^{\frac{2\alpha}{\alpha+1}}} P(y,s),
\ee
where $y=x/ (-t)^{\frac{1}{\a+1}}$ and $s=-\log (-t)$.
then the discrete self-similarity of the solution $(v,p)$ is equivalent to that $(V,P)$ is a solution to
\begin{equation}\label{se} \left\{
\aligned
&V_s+ \frac {\a}{\a+1} V +\frac {1}{\a+1}(y \cdot \nabla)V + (V\cdot \nabla )V =-\nabla
P,\\
& \mathrm{div}\,V=0,
\endaligned \right.
\end{equation}
which satisfies
$V(y,s)= V(y, s+S_0)$ for all $(y,s)\in \Bbb R^{3+1}$ with $S_0= -(\a+1)\log \lambda.$
The following result on the nonexistence of nontrivial discretely self-similar blow-up for the 3D Euler equations is proved in  \cite{cha-tsa}.
\begin{thm}
\label{th:1}
Let $V(y,s)\in C^1(\R^{3+1})$ be a time
periodic solution of \eqref{se} with period ${\sz}>0$.
\begin{itemize}
\item[(i)] Let $V$ satisfies either one of the following conditions:
\begin{itemize}
\item[(a)]  $V(y,s)\in C^2_y C^1_s (\Bbb R^{3+1}), $ and
$$ \Omega:=\mathrm{curl} \,
V\in \cup_{r>0}\cap _{0<q<r} L^q (\Bbb R^{3}\times [0, S_0]).
$$
\item[(b)]$\lim_{|y|\to \infty} \sup_{0 <s <\sz} | \nabla V(y,s)| =0,$ and $ \O\in L^q (\R^3\times [0,\sz])
    $ for some $q\in (0, \frac{3}{1+\a})$.
    \end{itemize}
    Then, $\O=0$. If we further assume $\quad \lim_{|y|\to \infty} V(y,s)=0\, \forall s\in [0, S_0)$, then $V=0$.
    \item[(ii)]
 Let $-1< \alpha \le 3/p$ or $3/2< \alpha < \I$, and  $V \in
L^p(\R^3 \times ([0,\sz])$ for some $3 \le p \le \I$. We also  assume that the pressure is given by $-\Delta_y P(\cdot,s) = \sum_{i,j} \pd_i \pd_j (V_iV_j(\cdot,s))$. Then, $V=0$.
\end{itemize}
\end{thm}
In the next section we will prove following theorem, which is an extension of the part (i)(b) of  the above theorem.
\begin{thm}
 Let $V(y,s)\in C^1(\R^{3+1})$ be a time
periodic solution of \eqref{se} with period ${\sz}>0$  satisfying
\bb\label{decayassume}\lim_{|y|\to \infty} \sup_{0 <s <\sz} | \nabla V(y,s)|=0,
\ee
and
\bb\label{th1.4} \int_0 ^{S_0} \int_{\Bbb R^3} |\O (y,s)|^q (1+|y|)^\eta  dyds <\infty, 
\ee
where $q\in (0, \frac{3+\eta}{1+\alpha})$ with $\alpha >-1, \eta >-3$.  Then, $\O=0$. Therefore if we further assume$\quad \lim_{|y|\to \infty} V(y,s)=0\, \forall s\in [0, S_0)$, then $V=0$.
\end{thm}
Note that case of $\eta=0$  of the above theorem corresponds to the part (i)(b) of Theorem 1.3.  We now  introduce more general notion of the asymptotically discretely self-similar blow-up for the 3D incompressible Euler equations.
Let $v\in C^1(\Bbb R^{3+1})$ be a solution to (E), which blows up at $t=0$.
We  say  that $t=0$ is a time for  {\em asymptotically discretely self-similar blow-up }  if
there exist  $\a >-1$ and a solenoidal vector field  $\overline{V }\in C^1(\Bbb R^{3+1} )$, which is a time periodic solution to (\ref{se}), and  satisfies the following convergence conditions:
\bb\label{conv2}
\lim_{t\uparrow 0}(-t)\left\| \nabla v(\cdot ,t) +\frac{1}{t} \nabla \overline{V}\left(\frac{\cdot}{(-t)^{\frac{1}{\alpha+1}}}, -\log (-t)\right)\right\|_{L^\infty }=0,
\ee
and there exists $\vare_0>0$ such that
\bb\label{conv3}
\sup_{-\vare_0 <t <0}(-t)^{{\frac{3-\a}{\a +1}}}\left\|  v(\cdot ,t) -\frac{1}{(-t)^\frac{\a}{\a+1}} \overline{V}\left(\frac{\cdot}{(-t)^{\frac{1}{\alpha+1}}}, -\log (-t)\right)\right\|_{L^1}<\infty.
\ee
The following is a generalization of the corresponding  results in \cite{cha1, cha2a}.
\begin{thm}
Let $T>0$ and  $v\in C([-T, 0); H^m (\Bbb R^3)) $, $m>5/2$,  be a classical solution to (E), which blows up at $t=0$.
Then, $t=0$ is not a time for  asymptotically discretely self-similar blow-up if  the blow-up profile $\overline{ V}\in   C^1(\Bbb R^{3+1})$ satisfies  either the conditions  (i) (a) or (i)(b) of Theorem 1.3,  or the conditions (\ref{decayassume})-(\ref{th1.4}) of Theorem 1.4, 
together with  the condition $\lim_{|y|\to \infty} \sup_{s\in \Bbb R} |V(y,s)|=0$.
\end{thm}

\section{Proof of the main theorems}

\noindent{\bf Proof of Theorem 1.1 } We recall the main  result in \cite{ser}, saying that $t=0$ is the blow-up time for the solution $v\in C(-\infty, 0; L^3 (\Bbb R^3))$  only if
\bb\label{sere}\lim_{t\uparrow 0} \|v(\cdot,t)\|_{L^3}=\infty.
\ee
 By the discrete self-similarity there exists $0<\lambda\neq 1$ such that
$$\lambda^k v(\lambda ^k x, \lambda ^{2k} t)= v(x,t)\quad \forall k\in \Bbb Z. $$
Hence,
\bb \|v(\cdot, \lambda ^{2k} t)\|_{L^3}= \|v(\cdot, t)\|_{L^3} \quad \forall t\in (-\infty, 0).\ee
Passing $k\to \infty $ if $\lambda \in (0, 1)$, while passing $k\to -\infty$ if $\lambda \in (1, \infty)$, one has for $t_k =\lambda^k t$,
$$\|v(\cdot ,t)\|_{L^3}=\lim_{t_k\to 0 } \|v(\cdot, t_k)\|_{L^3} =
\infty\quad \forall t\in (-\infty, 0), $$
 which is  contradiction to the fact that
 $v\in C(-\infty, 0; L^3 (\Bbb R^3))$. $\square$\\
\ \\
For the proof of Theorem 1.2 we  shall use the following result, which is Theorem 1.1 of \cite{gus}.
 \begin{lem} Let $q\in [3/2, \infty]$. Suppose $v$ is a suitable
 weak solution of (NS) in a cylinder, say  $Q=B(0,r_1)\times (-r_1^2, 0)$ for some $r_1>0$.
 Then, there exists a constant
 $\eta=\eta(q) >0$ such that if
 \bb\label{kang}
 \lim\sup_{r\downarrow 0}\left\{ r^{\frac{q-3}{q}}\mathrm{ess}\sup_{-r^2<t<0}\|v(\cdot,t)\|_{L^q (B(0, r))} \right\}\leq \eta,
 \ee
 then   $v$ is H\"{o}lder continuous both in space and time
 variables near $(0,0)\in \Bbb R^{3+1} $.
 \end{lem}

\noindent{\bf Proof of Theorem 1.2 }
We use the self-similar variables
$ y= x/\sqrt{-t}, \quad s=-\log(-t),$ and transform $(v,p)\to (V,P)$ as previously,
\bb
 v(x,t)= \frac{1}{\sqrt{-t}} V\left(\frac{x}{\sqrt{-t}}, -\log(-t)\right),\quad
 p(x,t)= -\frac{1}{t} P\left(\frac{x}{\sqrt{-t}}, -\log(-t)\right).
\ee
Then, substituting $(v,p)$ into (NS), we find that $(V,P)$ satisfies (\ref{leray}).
 We observe that the condition (\ref{conv1}) for some $R\in
(0, \infty)$ is equivalent to
 \bq
\label{convv14}
  \lim_{t\uparrow 0}(-t)^{\frac{q-3}{2q}}\sup_{t<\tau<0}\left\|v(\cdot, \tau)
-\frac{1}{\sqrt{-\tau}} \bar{V} \left(\frac{\cdot }{\sqrt{-\tau}}, -\log(-\tau)
 \right)\right\|_{L^q(B(0,R\sqrt{-t}\,))}
=0
  \eq
{\em for all } $R\in (0, \infty)$(see e.g. the argument in the proof of Theorem 1.5, pp. 446, \cite{cha1}), which implies,  in terms of $V$, that
 \bb\label{convv15}
 \lim_{s\to \infty} \|V(\cdot, s)-\bar{V}(\cdot , s)\|_{L^q(B(0, R))}=0
 \ee
 for all $R\in (0, \infty)$.
  We now show that this implies that $\overline{V}\in  C^1_tC^2_x (\Bbb R^{3+1})$ is a solution to (\ref{leray}).
Let $\xi
\in C^1_0 (0,S_0)$  and $\phi =(\phi_1 , \phi_2, \phi_3 )\in C_0^2 (\Bbb
 R^3)$ with div $\phi=0$. We multiply the first  equations of (\ref{leray}) by $\xi
 (s-S_0n)\phi (y)$, and integrate it over $\Bbb R^3\times [n, n+S_0]$,
 and then we  integrate by parts to obtain:
 \bq\label{wk1}
 &&-\int_0^{S_0}\int_{\Bbb R^3} \xi _s(s) \phi(y)\cdot V(y,s+S_0n)dyds
-\int_0 ^{S_0}\int_{\Bbb R^3}\xi (s)  V(y, s+S_0n)\cdot\phi(y) dyds \n \\
 &&\qquad-\frac{1}{2}\int_0 ^{S_0}\int_{\Bbb R^3}\xi (s)V(y, s+S_0n)\cdot (y \cdot \nabla)\phi (y)dyds\n \\
 &&\qquad -\int_0 ^{S_0}\int_{\Bbb R^3}\xi (s)\left[V(y,s+S_0n)\cdot (V(y,s+S_0n)\cdot \nabla )\phi (y)\right]  dyds
 \n  \\
 &&\quad= \int_0 ^{S_0}\int_{\Bbb R^3}\xi (s)  V(y, s+S_0n)\cdot\Delta \phi(y) dyds
 \eq
Similarly we multiply the second equations of (\ref{leray})  by $\xi
 (s-S_0n)\psi (y)$ with $\psi \in C^1_0 (\Bbb R^3)$, and integrate it over $\Bbb R^3\times [n, n+S_0]$, we have
 \bb\label{wk2}
 \int_0 ^{S_o} \int_{\Bbb R^3} V(s+S_0n)\cdot \nabla \psi (y)\xi(s) dyds=0.
 \ee
  Passing to the limit $n\to \infty$ in (\ref{wk1})-(\ref{wk2}), and  observing that $V(\cdot, s+S_0n)\to \bar{V}(\cdot, s)$ in $L^q_{loc} (\Bbb R^3) \hookrightarrow L^2_{\mathrm{loc}} (\Bbb R^3)$ for $q\in [2, \infty]$,
  we  find that $\bar{V}\in  C^1_sC^2_y (\Bbb R^{3+1}))$ satisfies (\ref{wk1}) and (\ref{wk2}).
  Integrating by part in  (\ref{wk1}) and (\ref{wk2}) for $\overline{ V}$  respectively, we obtain
\bb
\int_{\Bbb R^3}\int_0 ^{S_0} \left[ \overline{V_s}+ \frac {1}{2} \overline{V} +\frac {1}{2}(y \cdot \nabla)\overline{V} + (\overline{V}\cdot \nabla )\overline{V}-\Delta \overline{V}\right]\cdot \phi(y)\xi(s)dyds=0
\ee
for all vector test function $\phi \in C^2_0(\Bbb R^3)$ with div $\phi=0$, and $\xi \in C^1_0 (0, S_0)$, and also for all $\psi \in C^2(\Bbb R^3)$ we have
\bb
 \int_0 ^{S_0}\int_{\Bbb R^3} [\mathrm{div} \,\bar{V}]\psi (y)\xi(s)dyds=0.
 \ee
 Therefore there exists $\bar{P}\in C^1 (\Bbb R^3\times [0, S_0))$ such that
 $$
\overline{V}_s+ \frac {1}{2} \overline{V} +\frac {1}{2}(y \cdot \nabla)\overline{V} + (\overline{V}\cdot \nabla )\overline{V}-\Delta \overline{V}=-\nabla \overline{P},\quad
  \mathrm{div} \, \overline{V}=0.
$$
Since $\overline{V}$ is a $ C(\Bbb R^3; L^3(\Bbb R^3)\cap C^2 (\Bbb R^3))$ solution of (\ref{leray}) satisfying the periodic condition, $\overline{V}(y,s)=\overline{V}(y, s+S_0)$,  we find that $\overline{V}=0$ by Theorem 1.1.
Therefore, the assumption (\ref{conv1}) is reduced to
 \bb \label{pro15}
  \lim_{t\uparrow 0}\left\{(-t)^{\frac{q-3}{2q}}\sup_{t<\tau<0}\left\|v(\cdot,
\tau)\right\|_{L^q(B(3,R\sqrt{-t}\,))}\right\} =0
  \ee
  for all $R\in (0, \infty)$.
 We set $R=1$ and $\sqrt{-t}=r$ in
  (\ref{pro15}), then we obtain
 \bb \label{pro16}
\lim_{r\downarrow 0}\left\{r^{\frac{q-3}{q}}\sup_{-r^2<\tau<0}\left\|v(\cdot,
\tau )\right\|_{L^q(B(0,r))}\right\} =0.
  \ee
Applying Lemma 2.1,  we find that the space-time point $(0,0)$ is not a blow-up space-time point. $\square$\\
\ \\
In order to prove Theorem 1.3 we recall the following blow-up criterion for the Euler equations, which corresponds to Lemma 2.1 for the Navier-Stokes equations.
\begin{lem}
Let $m>5/2$, $\vare >0$ and $v\in C((-\vare,0); H^m (\Bbb R^3))$ be a classical solution to (E). Suppose the following inequality holds
\bb
\lim\sup_{t\uparrow 0}(-t)\|\nabla v(t)\|_{L^\infty}<1,
\ee
then there exists no blow-up at $t=0$.
\end{lem}
For the proof see the proof of Theorem 1.1 in \cite{cha2a}(see also \cite{hou} for similar result with different approach).\\

\noindent{\bf Proof of Theorem 1.5 }
We transform the solution $(v,p)\to (V,P)$ as in (\ref{ss}). Then, $(V,P)$ satisfies (\ref{se}). The convergence conditions (\ref{conv2})-(\ref{conv3}) can be written in  the self-similar form as
\bb\label{conva}
\lim_{s\to \infty} \|\nabla V(\cdot, s)-\nabla \overline{V}(\cdot, s)\|_{L^\infty}=0,
\ee
and
\bb\label{convb}
\sup_{-\log(\vare_0)<s<\infty}\| V(\cdot, s)- \overline{V}(\cdot, s)\|_{L^1}<\infty
\ee
respectively. From (\ref{conva}) and (\ref{convb}), using the interpolation for the $L^p$ spaces,  we have that
\bb
\lim_{s\to \infty} \|V(\cdot, s)- \overline{V}(\cdot, s)\|_{L^2 (B_R)}=0\quad  \forall R>0,
\ee
and repeating the argument of the proof of Theorem 1.2 word by word, we find that $\overline{V}$ is $ C^1(\Bbb R^3\times (0, S_0))$ solution of (\ref{se}), satisfying the time periodicity and one of the conditions  (i) (a) or (i)(b) of Theorem 1.3,  or the conditions (\ref{decayassume})-(\ref{th1.4}) of Theorem 1.4, 
together with  the condition $\lim_{|y|\to \infty} \sup_{s\in \Bbb R} |V(y,s)|=0$. 
Therefore  $\overline{V}=0$, and the condition (\ref{conv2}) reduces to
$$ \lim_{t\uparrow 0} (-t) \|\nabla v (\cdot, t)\|_{L^\infty}=0. $$
Thanks Lemma 2.2 we can conclude that $t=0$ is not a blow-up time. $\square$\\

\noindent{\bf Proof of Theorem 1.4 } We consider the vorticity equation of (\ref{se}),
\begin{equation}
\label{vse} \left\{
\aligned
&\O_s+  \O +\frac {1}{\a+1}(y \cdot \nabla)\O + (V\cdot \nabla )\O -(\O\cdot \nabla) V =0,\\
& \mathrm{div}\,V=0,\quad \mathrm{curl }\, V=\O.
\endaligned \right.
\end{equation}
We introduce a cut-off function $\sigma\in C_0
^\infty(\Bbb R^N)$ such that
\begin{equation}\label{16}
   \sigma(|x|)=\left\{ \aligned
                  &1 \quad\mbox{if $|x|<1$},\\
                     &0 \quad\mbox{if $|x|>2$},
                      \endaligned \right.
\end{equation}
and $0\leq \sigma (x)\leq 1$ for $1<|x|<2$.  For each $R>0$, we
define $\s_R (x):= \s \left(\frac{|x|}{R}\right)$.
Given $\rho>0$, we also define a function $\psi=\psi_\rho (y)$  as follows.
$$
\psi_\rho (y)=\left\{ \aligned &1,  &\mbox{if}\quad |y|>\rho+\pi\\
&\frac12\sin\left(|y|-\rho-\frac{\pi}{2}\right)+\frac12, &\mbox{if}\quad\rho<|y|\leq \rho+\pi\\
                                &0,        &\mbox{if} \quad|y|\leq \rho.
                                \endaligned
                                \right.
                                   $$
      We take $L^2 (\Bbb R^3\times[0, S_0])$ inner product the first equations of (\ref{vse}) by $\O |\O|^{q-2} |y|^\eta \psi_\rho \si$, and integration by parts to obtain
 \bq\label{In}
\lefteqn{0=\left(1-\frac{3+\eta}{q(\a +1)}\right)\int_0 ^{S_0} \int_{\Bbb R^3} |\O|^q  |y|^\eta \psi_\rho \si dyds
-\frac{1}{\a+1} \int_0 ^{S_0} \int_{\Bbb R^3} |\O|^q  |y|^\eta (y\cdot \nabla )\psi_\rho \si dyds}\n \\
&&-\frac{1}{\a+1} \int_0 ^{S_0} \int_{\Bbb R^3} |\O|^q  |y|^\eta\psi_\rho (y\cdot \nabla )\si dyds
-\frac{\eta}{q} \int_0 ^{S_0} \int_{\Bbb R^3} |\O|^q   V\cdot \frac{y}{|y|} |y|^{\eta-1} \psi_\rho \si dyds \n \\
&&-\frac{1}{q} \int_0 ^{S_0} \int_{\Bbb R^3} |\O|^q  |y|^\eta(V\cdot \nabla )\psi_\rho \si dyds
-\frac{1}{q} \int_0 ^{S_0} \int_{\Bbb R^3} |\O|^q  |y|^\eta\psi_\rho (V\cdot \nabla )\si dyds \n \\
&&\qquad-  \int_0 ^{S_0} \int_{\Bbb R^3} (\O\cdot \nabla) V\cdot \O |\O|^{q-2} |y|^\eta\psi_\rho \si dyds    \n \\
&&:=I_1+\cdots+I_7 ,
\eq
where the  periodicity of $V(y,s)$ and $\O(y,s)$ in $s$ was used.
We first observe that for all $\eta \in \Bbb R$
\bb\label{domi}
\int_0 ^{S_0} \int_{\Bbb R^3} |\O|^q  |y|^\eta \psi_\rho dyds
\leq \int_0 ^{S_0} \int_{\Bbb R^3} |\O|^q (1+ |y|)^\eta  dyds<\infty.
\ee
Therefore,  we can apply the dominated convergence theorem in the following.
\bq
|I_3|&\leq& \frac{\|\nabla \s\|_{L^\infty}}{R(\a+1)} \int_0 ^{S_0} \int_{\{ R\leq |y|\leq 2R\}} |\O|^q  |y|^\eta \psi_\rho  |y| dyds\n \\
&\leq& \frac{2\|\nabla \s\|_{L^\infty}}{(\a+1)} \int_0 ^{S_0} \int_{\{ R\leq |y|\leq 2R\}} |\O|^q  |y|^\eta \psi_\rho dyds\n \\
&\to& 0
\eq
as $R\to \infty$, and
\bq
|I_6|&\leq& \frac{\|\nabla \s\|_{L^\infty}}{q R} \int_0 ^{S_0} \int_{\{ R\leq |y|\leq 2R\}} |\O|^q  |y|^\eta \psi_\rho |V| dyds\n \\
&\leq& \frac{\|\nabla \s\|_{L^\infty}}{q R} \sup_{R\leq |y|\leq 2R} \sup_{s\in [0, S_0]}\frac{|V(y,s)|}{|y|}\int_0 ^{S_0} \int_{\{ R\leq |y|\leq 2R\}} |\O|^q  |y|^\eta \psi_\rho |y| dyds\n \\
&\leq& \frac{2\|\nabla \s\|_{L^\infty}}{q } \sup_{R\leq |y|\leq 2R} \sup_{s\in [0, S_0]}\frac{|V(y,s)|}{|y|}\int_0 ^{S_0} \int_{\{ R\leq |y|\leq 2R\}} |\O|^q  |y|^\eta \psi_\rho  dyds\n \\
&\to& 0
\eq
as $R\to \infty$. Here we used the following simple observation: Let $0\not= y_1\in \Bbb R^3$. Then,
$$V(y,s)=V(y_1,s) +\int_0 ^1 \frac{d}{d\tau} V(\tau y+(1-\tau )y_1 ) d\tau= V(y_1,s) +\int_0 ^1 y\cdot \nabla  V(\tau y+(1-\tau )y_1 ) d\tau,$$
and therefore,
$$
\frac{|V(y,s)|}{|y|}\leq \frac{|V(y_1,s)|}{|y|}+\sup_{0<\tau <1}  |\nabla V(\tau y+(1-\tau )y_1 ,s)|
 \to 0\quad \mbox{as}\quad |y|\to \infty,
$$
which follows from the hypothesis (\ref{decayassume}). Therefore, passing $R\to \infty$, and applying the dominated convergence theorem to the other terms of (\ref{In}), we find that
\bq\label{Ina}
\lefteqn{0=\left(1-\frac{3+\eta}{q(\a +1)}\right)\int_0 ^{S_0} \int_{\Bbb R^3} |\O|^q  |y|^\eta \psi_\rho  dyds
-\frac{1}{\a+1} \int_0 ^{S_0} \int_{\Bbb R^3} |\O|^q  |y|^\eta (y\cdot \nabla )\psi_\rho  dyds}\n \\
&&-\frac{\eta}{q} \int_0 ^{S_0} \int_{\Bbb R^3} |\O|^q   V\cdot \frac{y}{|y|} |y|^{\eta-1} \psi_\rho dyds
-\frac{1}{q} \int_0 ^{S_0} \int_{\Bbb R^3} |\O|^q  |y|^\eta(V\cdot \nabla )\psi_\rho  dyds\n \\
&&\qquad-  \int_0 ^{S_0} \int_{\Bbb R^3} (\O\cdot \nabla) V\cdot \O |\O|^{q-2} |y|^\eta\psi_\rho  dyds    \n \\
&&:=J_1+\cdots+J_5.
\eq
Under our hypothesis we have $J_1 \leq 0$. Since $\psi_\rho(y)$ is radially non-decreasing, we also have $J_2 \leq0$, and
\bb\label{j2}
J_2=-\frac{1}{\a+1} \int_0 ^{S_0} \int_{\{\rho \leq |y|\leq \rho+\pi\}} |\O|^q  |y|^{\eta+1} \psi'_\rho  dyds.
\ee
 We have
\bq\label{j3}
|J_3|\leq \frac{|\eta|}{q}\sup_{|y|\geq \rho, s\in [0, S_0]} \frac{|V(y,s)|}{|y|} \int_0 ^{S_0} \int_{\Bbb R^3} |\O|^q    |y|^{\eta} \psi_\rho dyds\leq& \frac14 |J_1|
\eq
for sufficiently large $\rho$. We compute
\bq\label{j4}
|J_4|&\leq & \frac{1}{q} \int_0 ^{S_0} \int_{\Bbb R^3} |\O|^q  |y|^\eta|V|\psi'_\rho  dyds\n \\
&\leq& \frac{1}{q}\sup_{|y|\geq \rho, s\in [0, S_0]} \frac{|V(y,s)|}{|y|}\int_0 ^{S_0} \int_{\{\rho \leq |y|\leq \rho+\pi\}} |\O|^q|y|^{\eta+1} \psi'_\rho  dyds\n \\
&\leq& \frac12 |J_2|
\eq
for sufficiently large $\rho$. For $J_5$ we obtain
\bq\label{j5}
|J_5|&\leq& \int_0 ^{S_0} \int_{\Bbb R^3} |\nabla V| |\O|^q |y|^\eta\psi_\rho  dyds    \n \\
&\leq& \sup_{|y|\geq \rho, s\in [0, S_0]} |\nabla V(y,s)| \int_0 ^{S_0} \int_{\Bbb R^3}|\O|^q |y|^\eta\psi_\rho  dyds\n \\
&\leq & \frac14 |J_1|
\eq
for sufficiently large $\rho$. Taking into account the estimates (\ref{j3})-(\ref{j5}) in (\ref{Ina}), we find that there exists $\rho_0>0$
such that
\bq \lefteqn{0\geq  \frac12 \left(\frac{3+\eta}{q(\a +1)}-1\right)\int_0 ^{S_0} \int_{\Bbb R^3} |\O|^q  |y|^\eta \psi_\rho  dyds}\n \\
&&\qquad+\frac{1}{2(\a+1)} \int_0 ^{S_0} \int_{\Bbb R^3} |\O|^q  |y|^\eta (y\cdot \nabla )\psi_\rho  dyds
\eq
for all $\rho \geq \rho_0$. Hence, $\O(y,s)=0$ for all $(y,s)\in \{ y\in \Bbb R^3\, |\, |y|>\rho_0\}\times [0, S_0]$.
Applying Theorem 1.3 (i)(a), we conclude $\O=0$. $\square$\\
$$ {\bf Acknowledgements} $$
The author would like to thank  to the anonymous referees for useful suggestions.
  This research is supported partially by NRF Grant 2006-0093854 and 2009-0083521.
  
\end{document}